\documentclass[12pt,a4paper,leqno]{article}

\usepackage{fullpage,amsmath,amsthm,amsfonts,amscd,graphicx,latexsym,amssymb,eucal,curves}
\usepackage[all]{xy}

\newcommand {\const}{\text{\rm const.}\;}

\newcommand {\Pfeil} {\longrightarrow}

\newcommand{\An}{\underline{\rm An}}

\newcommand{\abs}[1]{\lvert#1\rvert}

\newcommand {\und}[1]{\underline#1}

\renewcommand{\phi}{\varphi}

\renewcommand{\epsilon}{\varepsilon}

\newcommand{\BS}{{\rm BS}}
\newcommand{\ab}{{\rm ab}}
\newcommand{\alg}{{\rm alg}}
\newcommand{\an}{{\rm an}}

\newcommand{\fin}{{\rm fin}}

\newcommand{\orb}{{\rm orb}}

\newcommand{\sing}{{\rm sing}}
\newcommand{\temp}{{\rm temp}}
\renewcommand{\top}{{\rm top}}

\newcommand{\cov}{{\rm cov}}

\newcommand{\Set}{{\rm \underline{Sets}}}

\DeclareMathOperator{\im}{im}
\DeclareMathOperator{\Aut}{Aut}
\DeclareMathOperator{\Out}{Out}

\DeclareMathOperator{\Cov}{Cov}

\DeclareMathOperator{\MOrb}{MOrb}

\DeclareMathOperator{\Inn}{Inn}

\DeclareMathOperator{\PGL}{PGL}

\DeclareMathOperator{\rk}{rk}

\DeclareMathOperator{\val}{val}

\newcommand{\porbS}{{\pi_1^\orb(\mathcal{S},\bar{s})}}
\newcommand{\ptempS}{{\pi_1^\temp(S,\bar{s})}}
\newcommand{\ptopS}{{\pi_1^\top(S,\bar{s})}}

\newcommand{\Res}{{\rm Res}}

\newtheorem{Def}{Definition}[section]

\newtheorem{Rem}[Def]{Remark}

\newtheorem{Prop}[Def]{Proposition}

\newtheorem*{Satz*}{Satz}

\newtheorem{Lemma}[Def]{Lemma}

\newtheorem{Cor}[Def]{Corollary}

\newtheorem*{Quest*}{Question}

\newtheorem{Thm}{Theorem}

\newcommand{\Z}{\mathbb{Z}}

\newcommand{\C}{\mathbb{C}}

\newcommand{\F}{\mathbb{F}}
\newcommand{\G}{\mathbb{G}}

\renewcommand{\P}{\mathbb{P}}

\newcommand {\fB} {\mathfrak B}

\newcommand {\fH} {\mathfrak H}

\newcommand {\fK} {\mathfrak K}

\newcommand {\fM} {\mathfrak M}

\newcommand {\fS} {\mathfrak S}
\newcommand {\fT} {\mathfrak T}
\newcommand {\fU} {\mathfrak U}

\newcommand {\cC} {\mathcal C}

\newcommand {\cM} {\mathcal M}

\newcommand {\cS} {\mathcal S}

\newcommand {\CC} {\mathbb C}

\newcommand {\Cp}{\CC_p}

\begin{document}



\title{\bf Construction of $p$-adic Hurwitz spaces}

\author{Patrick Erik Bradley}
\date{}

\maketitle






\begin{abstract}
Moduli spaces for Galois covers of $p$-adic Mumford curves by 
Mumford curves are constructed using Herrlich's Teichm\"uller
spaces, Andr\'e's orbifold fundamental groups, and Kato's graphs
of groups encoding ramification data of charts for Mumford orbifolds.

\end{abstract}

\section{Introduction}

Hurwitz spaces are moduli spaces for finite branched covers of curves.
Subspaces of their $p$-adic versions are moduli spaces for 
covers of
{\em Mumford curves},
i.e.~Schottky uniformisable curves.
In this paper we give a description of special kinds of such 
spaces parametrising Galois covers of Mumford curves by Mumford curves.

\medskip
These so-called {\em Mumford-Hurwitz spaces} turn out to be a
finite disjoint union of equidimensional
 moduli spaces of $N$-uniformisable Mumford curves,
where $N$ runs through some special type of
finitely generated groups: Bass-Serre fundamental groups of
graphs of groups. For such moduli spaces, there  exists a theory of
Teichm\"uller spaces \cite{Her85}.
This  allows to ``bundle together'' those moduli spaces involved
by examining which quotients of  Y.~Andr\'e's 
 orbifold
fundamental group 
are Bass-Serre fundamental groups for charts of a given reference 
orbifold. The finiteness of the number of components then boils
down to giving bounds for the number of vertices of possible
graphs.  

\medskip
Equidimensionality of all components is achieved by relating the number
of branch points to the dimension of the Teichm\"uller spaces involved.
This $p$-adic counterpart of the well-known
complex result does not seem to have been mentioned in the literature, so far.

\medskip
We adopt Berkovich's language of $p$-adic strictly analytic spaces
\cite{Ber90} and make use of Andr\'e's theory
of orbifold fundamental groups from \cite{AndIII}
which was developed for studying differential equations
on one-dimensional orbifolds allowing covers by Mumford
curves: we call these, following \cite{Kat01}, {\em Mumford orbifolds}.

\medskip
In Section 2, we give a very brief overview of the theory of temperate
fundamental groups, $p$-adic orbifolds and their fundamental groups.

\smallskip
Section 3 reviews Teichm\"uller spaces and moduli spaces of Mumford
curves
from the point of view of $p$-adic geometry. It reveals that
Herrlich's Teichm\"uller spaces parametrise uniformisations of 
Mumford orbifolds.

\smallskip
In Section 4, the actual construction of Mumford-Hurwitz spaces
is performed, and the finiteness result is proven: a main ingredient
is the quotient pasting of Kato trees, used in \cite{Brad02a}
for obtaining formulae for the number of branch points.
The section ends with some examples.

\medskip
The results of this paper can be found in the author's dissertation
\cite{Brad02} of which
we recommend sections 2 and 3 as a quick introduction to the theory
of $p$-adic analytic spaces and to the temperate
and orbifold fundamental groups from \cite{AndIII}.

\section{Orbifolds and orbifold fundamental groups}

\subsection{The temperate fundamental group}

Following \cite{AndIII}, we have

\begin{Def}
Let $S$ be a $p$-adic manifold. A {\em geometric point} of $S$ is an
analytic
map $\bar{s}\colon M(\Omega)\to S$, where $(\Omega,|\; |)$ is
an algebraically closed, complete extension of $(\Cp,|\; |_p)$.
\end{Def}
As $\Omega$ is a normed $\Cp$-algebra, $M(\Omega)$ is well-defined.
In fact, it is a point, and with $s:=\bar{s}(M(\Omega))\in S$
we shall always mean the image of the geometric point $\bar{s}$.

\bigskip
In $p$-adic geometry, there are not enough topological covers of a
given manifold. A remedy for this, is to extend the definition
of ``cover'' \cite{dJ95}: a {\em cover} of $S$ is an analytic map
$f\colon S'\to S$, such that there is an open covering $\fU$ of $S$
with the property
$$
\forall\, U\in\fU\colon f^{-1}(U)=\coprod\limits_i V_i
$$
and all restrictions $f|_{V_i}\colon V_i\to U$ are finite.
We see that in case all $f|_{V_i}$ are isomorphisms, we recover the 
well-known notion of ``topological cover''.

\smallskip
A {\em morphism} of covers of $S$ is a commuting triangle
$$
\xymatrix{
S' \ar[rr]^\an \ar[dr]_\cov && S'' \ar[dl]^\cov\\
&S&
}
$$
The resulting category is called $\Cov_S$.

\medskip
This is all a bit too general, but de Jong invents ``\'etale''
covers as a special case in order to obtain a first sensible
fundamental group. This turns out to be too big, so Andr\'e
becomes even more special:

\begin{Def}
A cover $g\colon S'\to S$ is called {\em temperate}, if there is a topological
cover $T'\to T$ and a commuting diagramme
$$
\xymatrix{
T' \ar[r]^{/R} \ar[d]_\top & S' \ar[d]^g\\
T \ar[r]_\alg &S
}
$$ 
where the upper horizontal arrow is a quotient over $S$, and the
lower horizontal arrow is a {\em finite \'etale} or {\em algebraic}
cover: it is a finite morphism $f$ and the $f|_{V_i}$ as above are 
\'etale maps.

\smallskip
The sub-category of $\Cov_S$
of temperate covers is denoted by $\Cov^\temp_S$.
\end{Def}

As fundamental groups act as permutations of fibres of covers,
we define the {\em geometric fibre} of a geometric point
$\bar{s}\colon M(\Omega)\to S$ of $S$ for the cover
$f\colon S'\to S$ to be the set of all lifts of $\bar{s}$
to geometric points of $S'$:
$$
f^{-1}(\bar{s}):=\left\{\bar{s}'\colon M(\Omega)\to S'\;\left|\;
	\raisebox{10mm}{
	\xymatrix{
	&S'\ar[d]^f\\
	M(\Omega) \ar[ru]^{\bar{s}'}\ar[r]_{\bar{s}}&S
	}}
\quad \text{commutes}\right.\right\}\,.
$$

The {\em fibre functor} is the covariant functor
$$
F_{S,\bar{s}}\colon\Cov_S\to\Set,\;(f\colon S'\to S)\mapsto f^{-1}(\bar{s}).
$$ 

\begin{Def}
The {\em temperate fundamental group} of $S$ with base point 
$\bar{s}\colon M(\Omega)\to S$ is 
$$
\ptempS:=\Aut F^\temp_{S,\bar{s}},
$$
the automorphism group of the restricted fibre functor 
$F^\temp_{S,\bar{s}}:=F_{S,\bar{s}}|_{\Cov^\temp_S}$.
\end{Def}

\begin{Rem}
The temperate fundamental group is a pro-discrete topological group.
A basis of open neighbourhoods of $1$ is the system of normal closed
subgroups $H$ such that $\ptempS/H$ is the Galois group
of the universal topological cover of some algebraic Galois
cover of $S$.

\smallskip
Its pro-finite completion is the algebraic fundamental group
(in case $S$ is the analytification of an algebraic variety, the
latter group coincides with the usual algebraic fundamental group).
\end{Rem}

\begin{proof}
The first part is \cite[2.1.5.]{AndIII}. The second follows from
\cite[1.4.8]{AndIII}, as the algebraic covers of $S$
form a full sub-category of $\Cov^\temp_S$.
\end{proof}

\subsection{The orbifold fundamental group}

\begin{Def}
A {\em $p$-adic orbifold} $\cS=(\bar{S},(Z_i,e_i))$ is a 
$p$-adic manifold $\bar{S}$ together with a locally finite family
$(Z_i)$ of irreducible divisors and natural numbers $e_i>0$,
and which can be covered by {\em orbifold charts}, i.e.~by analytic maps 
$\phi\colon W\to
V\subseteq \bar{S}$ with the properties
\begin{itemize}
\item $W$ is a $p$-adic manifold.
\item $\phi\in\Cov_V$.
\item Outside $Z:=\bigcup\limits_i Z_i\cup \bar{S}_\sing$, $\phi$
is a temperate Galois cover, ramified above $Z_i$ with index $e_i$.
\end{itemize} 

A {\em morphism} $f\colon \cS'\to\cS$ of orbifolds is a morphism
$\bar{S}'\to \bar{S}$ such that $\bar{S}'$ is covered with orbifold
charts
$\phi'\colon W'\to V'$ and there exist orbifold charts 
$\phi\colon W\to V\subseteq \bar{S}$ for which the diagramme
$$
\xymatrix{
W' \ar[d] \ar[r]^{\phi'} & V' \ar[r] & \bar{S}' \ar[d]^f \\
W \ar[r]_\phi & V \ar[r] &S
}
$$
commutes.
\end{Def}

The category of $p$-adic manifolds will be viewed as a full sub-category
of the category of $p$-adic orbifolds.

\begin{Def}
A morphism $f\colon \cS'\to \cS$ is called an {\em orbifold cover}
if the underlying analytic map $f\colon \bar{S}'\to\bar{S}$
is in $\Cov_{\bar{S}}$, and for every orbifold chart 
$\phi'\colon W'\to V'\subseteq \bar{S}'$ the composition 
$f\circ\phi\colon W'\to f(V')$ is an orbifold chart of $\bar{S}$.
\end{Def}

This gives us a category $\Cov_{\cS}$. If both $\cS'$ and $\cS$ are
manifolds, then an orbifold cover is nothing but a 
temperate cover.

\smallskip
If an orbifold  has a global orbifold chart, it will be called
{\em uniformisable}. Let  $S:=\bar{S}\setminus Z$,
where $Z:=\bigcup\limits_i Z_i\cup \bar{S}_\sing$.

\begin{Rem}	\label{orbtemp}
Let $\cS$ be uniformisable and $(f\colon \cS'\to\cS)\in\Cov_{\cS}$.
Then the restriction $f^{-1}(S)\to S$ is a temperate cover.
\end{Rem}

\begin{proof}
\cite[4.4.5]{AndIII}
\end{proof}

Important is the following notion. Let a discrete
group $G$ act properly discontinuously on a connected orbifold 
$\cS=(\bar{S},(Z_i,e_i))$. Then the 
{\em orbifold quotient} is the
orbifold
$$
\cS/G:=(\bar{S}/G,(D\cdot G,e_D\cdot \abs{G_D})),
$$
where $D$ runs through the $Z_i$ and the prime divisors of $\bar{S}$
with non-trivial stabiliser $G_D$, and $e_D=1$ for $D\notin\{Z_i\}$.

\bigskip
Let $\cS$ be a uniformisable connected orbifold. Then Remark \ref{orbtemp}
gives us a restriction functor
$$
\Res\colon\Cov_{\cS}\to\Cov_S^\temp,
$$
which is faithful \cite[4.5.2]{AndIII}.

\begin{Def}
Let $\bar{s}$ be a geometric point of $S=\bar{S}\setminus Z$ for the
uniformisable orbifold $\cS$. Then 
$$
\porbS:=\Aut F_{\cS,\bar{s}},
$$ 
where $F_{\cS,\bar{s}}:=F^\temp_{S,\bar{s}}|_{\Cov_{\cS}}$,
is called the {\em orbifold fundamental group} of $\cS$ with base
point $\bar{s}$.
\end{Def}

The orbifold fundamental group is a pro-discrete topological group,
because of
\cite[1.4.7]{AndIII}, and there exists a fundamental exact sequence
for finite Galois orbifold covers:

\begin{Lemma}	\label{fundexact}
Let $f\colon \cS'\stackrel{/G}{\Pfeil}\cS$ be a finite orbifold quotient
cover. Then for every point $\bar{s}'$ of the geometric fibre 
$f^{-1}(\bar{s})$,
the sequence
$$
1\to\pi_1^\orb(\cS',\bar{s}')\to\porbS\to G\to 1
$$
is exact.
\end{Lemma}

\begin{proof}
Composing an orbifold cover of $\cS'$ with $f$ gives an orbifold cover of 
$\cS$,
so there is a natural morphism of topological groups
$$
\alpha\colon\pi_1^\orb(\cS',\bar{s}')\to\porbS.
$$
For the injectivity of $\alpha$, the proof of \cite[1.4.12(b)]{AndIII}
carries over: the system of all connected Galois orbifold covers
factorising through $f$ is cofinal in $\Cov_{\cS}$. This implies
that there is an open injective morphism from
$\pi_1^\orb(\cS',\bar{s}')
=\lim\limits_{\longleftarrow\atop X'\in\Cov_{\cS'}}\Aut(X'/\cS')$
into 
$$
\lim\limits_{\longleftarrow\atop X'\in\Cov_{\cS'}} \Aut(X'/\cS)
=\lim\limits_{\longleftarrow\atop X\in\Cov_{\cS}} \Aut(X/\cS)
=\porbS.
$$
The arrow $\beta\colon\porbS\to G$ exists and is surjective, because $G$
is the automorphism group of an orbifold cover of $\cS$
and as such $G$ is a quotient of $\porbS$. Here, $\ker\beta$
is the set of all $\gamma\in\porbS$ fixing $f^{-1}(\bar{s})$ pointwise.
But this equals $\im\alpha$, as each $\gamma \in\im\alpha$ lifts
to a fibre automorphism of $\bar{s}'$.
\end{proof}

\begin{Rem}
In the following, we will use this exact sequence in the special case 
that $\cS=\bar{S}$
is a manifold: then $\porbS=\ptempS$.

\smallskip
For example, the orbifold $\cS=(\P^1,(0,2),(1,2),(\infty,2),(\lambda,2))$
with $\lambda\notin\{0,1,\infty\}$ is uniformisable by an elliptic
curve $E\stackrel{/C_2}{\Pfeil}\cS$ given by the equation
$y^2=x(x-1)(x-\lambda)$ (we assume $p\neq 2$).
The exact sequence 
$$
1\to\pi_1^\temp(E,1)\porbS\to C_2\to 1
$$
shows that
$$
\porbS\cong\begin{cases} \hat{Z}^2\rtimes C_2, 
		& \text{if $|\lambda||\lambda-1|=1$ (good reduction)}\\
			(\Z\times\hat\Z)\rtimes C_2,
		& \text{if $|\lambda||\lambda-1|\neq 1$ (Tate curve)}
	\end{cases}
$$
The temperate fundamental group of the Tate curve can be calculated
from the sequence
$$
1\to\pi_1^\temp(\G_m,1)\to\pi_1^\temp(E,1)\to\pi_1^\top(E,1)\to 1
$$
coming from the topological universal cover $\G_m\to E$, which can be
seen as a quotient orbifold cover with $\pi_1^\top(E,1)$
acting freely on $\G_m$ (that is why the quotient is also a manifold)
{\rm \cite[2.3.2]{AndIII}}.
In Section \ref{Tate}, we will study the action of $C_2$ on the temperate fundamental group
of an elliptic curve.
\end{Rem}

\begin{Def}
A {\em Mumford orbifold} is a one-dimensional orbifold $\cS$ uniformisable
by
a global Galois orbifold chart $\bar{S}'\to \cS$ with a Mumford curve 
$\bar{S}'$. If the orbifold $(\P^1,(0,e_0),(1,e_1),(\infty,e_\infty))$
is a Mumford orbifold, it will be called a {\em Mumford-Schwarz orbifold}.
\end{Def}

In the sequel, we will be dealing mostly with Mumford orbifolds.



\section{Teichm\"uller spaces}
\subsection{Mumford curves of genus $g\ge 2$}

The following theorem, of which we will give a
 short proof from
\cite[Satz 3.12]{Brad02} using the Berkovich geometry of one-dimensional
analytic spaces,
 is 
 well-known:

\begin{Thm}	\label{topfree}
If $S$ is an irreducible non-singular projective algebraic curve
defined over a large enough complete non-archimedean field $K$,
then $\ptopS$ is a finitely generated free group.
\end{Thm}

\begin{proof}
According to \cite[4.3.2]{Ber90}, $S$ is a special quasi-polyhedron,  its 
skeleton $\Delta(S)$ is a subgraph of the intersection graph
$\Delta^\an(S)$ of its stable reduction (this is where ``large enough'' 
enters), and $\Delta(S)$ has the same
Betti number as $\Delta^\an(S)$ (at most $g$).
Let $\Omega\to S$ be the topological universal cover. As $\ptopS$
also acts on $\Delta(S)$ \cite[4.1.8]{Ber90},
we have a commuting diagramme
$$
\xymatrix{
\Delta(\Omega) \ar@{^{(}->} [r] \ar [d]_{/\ptopS} &\Omega \ar[d]^{/\ptopS} \\
\Delta(S) \ar@{^{(}->} [r] & S
}
$$
Because of the retraction map $\Omega\to\Delta(S)$, this
diagramme is Cartesian. Therefore the left arrow going down is universal
in the category of graph covers of $\Delta(S)$. This implies
$\ptopS$ is the fundamental group  $\pi_1(\Delta(S)$, a free group
of rank at most $g$.   
\end{proof}

Let $C$ be a Mumford curve of genus $g\ge 2$ defined over $\Cp$.
Its stable reduction is a curve defined over $\bar{\F}_p$
and whose intersection graph has genus $g$. Let
$\Gamma:=\pi_1^\top(C,\bar{s})$, according to Theorem
\ref{topfree} a free group of rank $g$, and
let $\Omega\stackrel{/\Gamma}{\Pfeil} C$
be the topological universal cover. $\Omega$ is an analytic
subspace of $\P^1_{\Cp}$, and $\Gamma$ acts discontinuously on it. 
This gives a faithful representation $F_g\to\PGL_2(\Cp)$ of the fundamental
group as a Schottky group \cite{GvP80}.

\medskip
The set $\fT_g$ of faithful discontinuous 
representations $F_g\to\PGL_2(\Cp)$ is
known to be an open analytic sub-domain of the affine $\Cp$-variety
$\fS_g\cong\PGL_2(\Cp)^g$ of all representations of $F_g$ into
$\PGL_2(\Cp)$
 \cite[Sect.~1]{Her84}. $\PGL_2(\Cp)$ acts on $\fT_g$
by conjugation, and 
$$
\bar\fT_g:=\PGL_2(\Cp)/\fT_g
$$
is the {\em Teichm\"uller space} for $F_g$.

\smallskip
According to \cite{Her84}, $\Gamma_g:=\Out F_g$ acts on the Teichm\"uller 
space, and the quotient $\fM_g=\bar\fT_g/\Gamma_g$ is the moduli
space of Mumford curves of genus $g$, which can be viewed as an
analytic subspace of the moduli space  $\cM_g$
of all irreducible projective curves of genus $g$ \cite{Lut83}.

\medskip
In order to prove connectedness and simple connectedness of Teichm\"uller
space, Gerritzen dissects $\bar{\fT}_g$ into inadmissible
open parts $\fB_g(\Gamma)$ depending only on the possible stable
reduction graphs of Mumford curves. 
This gives the connectedness result for $\fM_g$ which inadmissibly
locally looks like $\fB_g(\Gamma)/\Aut\Gamma$ \cite{Ger81}.

\smallskip
From the Berkovich geometric 
viewpoint, the parts $\fM_g(\Gamma):=\fB_g(\Gamma)/\Aut\Gamma$
are not disjoint---only their sets of classical points are.
In fact,  intersections $\fM_g(\Gamma)\cap \fM_g(\Gamma')$
consist of generic points of discs $D$
corresponding to families of Mumford curves parametrised by $D$
whose skeletons are either $\Gamma$ or $\Gamma'$.

\subsection{Discontinuously uniformisable Mumford curves}

Generalising the results of the preceding section, Herrlich
constructs Teichm\"uller spaces for finitely generated groups $N$
\cite{Her85}. In our language we ought to proceed thus:
let $\fT(N)$ be the functor
$$
\An(\Cp)\to\Set, \; 
S\mapsto \{\psi\colon N\to\Aut(\P^1_S)\mid \text{$\psi$ is discontinuous}\}, 
$$ 
where a {\em discontinuous} representation $\psi$ over $S$ is meant to be
injective and for each geometric point $\bar{s}\colon K\to S$
the induced representation $\psi_{\bar{s}}\colon N \to \PGL_2(K)$
is discontinuous.

\begin{Prop}	\label{moduli}
The functor $\fT(N)$ is representable, i.e.~a fine moduli space,
and the following quotients
\begin{enumerate}
\item $\overline{\fT(N)}:=\PGL_2\backslash\fT(N)$, 
\item $\fM(N):=\overline{\fT(N)}/\Out(N)$ and 
\item $\fM(N,F):=\overline{\fT(N)}/\Out_F(N)$
\end{enumerate}
are well-defined, if $N$ contains a free subgroup  $F$ of rank $\ge 2$
and of finite index in $N$. In this case, the first one is representable,
and the other two
 lead to coarse moduli spaces.
\end{Prop}

Here, $\Out_F(N)$ means $\Aut_F(N)/\Inn(N)$, 
and $\Aut_F(N):=\{\alpha\in\Aut(N)\mid\alpha(F)=F\}$.

\begin{proof}
The first statement is Folgerung ii) on p.~149 in \cite{Her87}.

\smallskip
1.~is Folgerung on p.~151 in \cite{Her87}.

\smallskip
2.~is \cite[Satz 2]{Her87}. 

\smallskip
3.~can be proven the same way as 2.~with some slight modifications.
\end{proof}

\begin{Rem}
1. The coarse moduli space $\fM(N)$ parametrises so-called {\em $N$-uniformisable}
covers $\phi\colon\Omega'\stackrel{/N}{\Pfeil} C$ of Mumford curves of genus $\rk N^\ab$
which are temperate outside the branch locus.
If $(\zeta_i,e_i)$ are the branch points in $C$ together with the orders
of ramification, then $\phi$ is a global chart of the orbifold
$\cC=(C,(\zeta_i,e_i))$. So, the moduli space actually  parametrises all 
$N$-uniformisable Mumford orbifolds.

\smallskip
2. $\fM(F,N)$ parametrises commuting diagrammes
$$
\xymatrix{
\Omega' \ar[rr]^{/F} \ar[drr]_{/N}&& C'\ar[d]^{\text{\rm finite}}  \\
	&&C
}
$$
and is finite over $\fM(N)$.
\end{Rem}


\begin{Prop}	\label{equidim}
If $N$ contains a free group of rank $\ge 2$, then $\overline{\fT(N)}$
is a simply connected analytic space. Its finitely many components
are (smooth) $p$-adic manifolds of dimension
$$
3g-3 + \underbrace{2(C-c)+3(D-d)}_{=:n},
$$
where $n$ depends only on $N$ and equals the number of branch points of the
orbifold covers $\Omega(N)\to\cC$ in $\fM(N)$.
\end{Prop}

\begin{proof}
The rigid analytic version of this proposition is essentially
Herrlich's habilitation thesis \cite{Her85} or \cite{Her87},
except for the last statement. 

\smallskip
The number $n$ does not depend on the particular graph of groups with
fundamental group isomorphic to $N$: 

First,  $n$ is the number of 
cusps in the Kato graph of a realisation of $N$ 
as the Galois group of an orbifold cover $\phi\colon\Omega\to\cC$ of a
Mumford orbifold $\cC$
\cite[Theorem 2]{Brad02a}.
These cusps correspond bijectively to the branch points of $\phi$.
As all vertex groups are finite, Khramtsov's characterisation of
finite graphs of groups with isomorphic fundamental groups applies
\cite{Khr91}: two Kato graphs with the same fundamental group
differ by a finite number of admissible edge contractions or ``slides''
of an edge $e$ along an element $g\in N$ stabilising one of $e$'s
extremities $v$. The latter means that the topological graph and all 
vertex and edge groups are the same, 
only  the embedding 
$\alpha\colon N_e\to N_v$ is replaced by  
$c_g\circ \alpha$, where $c_g$ is conjugation by $g\in N_v$.  
The formula in \cite[Theorem 3]{Brad02a} now shows that the number of cusps
are the same for all Kato graphs with the same fundamental group $N$. 
\end{proof}

\section{Components of Hurwitz spaces}

\subsection{$G$-covers}

Let $\cS=(\bar S, (\zeta_i, e_i))$ be a Mumford orbifold, and $N$
a virtually free quotient of $\porbS$ of rank $\ge 2$ such that there exists
a faithful discontinuous representation
$\tau\colon N\to\PGL_2(\Cp)$. Let $\Gamma^*(\tau(N))$ be the Kato
graph for the global orbifold chart $\Omega^*(\tau(N))\to\cS$
associated to $\tau$.

\begin{Prop}
$\Cov^\BS_{\Gamma^*(\tau(N))}$ is a full sub-category of $\Cov_{\cS}$.
\end{Prop}

\begin{proof}
To a cover of graphs of groups
$$
\xymatrix{
\Delta^*(\tau(N)) \ar[dr]_{/N} \ar[r]& (\Gamma, H_\bullet)\ar[d]^\cov\\
&\Gamma^*(\tau(N))
}
$$
corresponds a unique subgroup $H$ of $\tau(N)$ giving rise to a unique
orbifold cover
$$
\xymatrix{
\Omega \ar[r]^{/H}\ar[dr]_{/\tau(N)}&\cS'\ar[d]^\orb\\
&\cS
}
$$
where $\cS'=(\Omega^*/H,(\zeta_{ij},e_{ij})$, with 
$\zeta_{ij}$ above $\zeta_i$ 
and $e_{ij}$ the ramification index of $\zeta_{ij}$ $(i=0,\dots,n)$.
This gives a functor $\Cov^\BS_{\Gamma^*(\tau(N))}\to\Cov_{\cS}$
embedding the first category into the second as a full sub-category.
\end{proof}

\'Etale covering theory yields

\begin{Cor}
There exists a canonical homomorphism of topological groups
$\psi\colon\porbS\to\pi_1^\BS(\Gamma^*(\tau(N)))$ with dense image. As 
$\pi_1^\BS(\Gamma^*(\tau(N)))$ is discrete, $\psi$ is surjective.
\end{Cor}

We observe that all (finite) Galois global orbifold charts 
$\phi\colon \bar{S}'\stackrel{/G}{\Pfeil}\cS$
 with a Mumford curve $\bar{S}'$ lead to covers of graphs:
let $\Omega'\stackrel{\top}{\Pfeil}\bar{S}'$ be the universal topological
cover of $\bar{S}'$, and let $N_\phi:=\Aut(\Omega',\cS)$. Then
$$
\Omega'\stackrel{/N_\phi}{\Pfeil}\cS
$$
is a Galois orbifold cover, and the following diagramme commutes with 
exact rows \cite[4.5.8]{AndIII}:
$$
\xymatrix{
1\ar[r] & \pi_1^\temp(\bar{S}',\bar{s}') \ar[r] \ar[d] & \pi_1^\orb(\mathcal{S},\bar{s}) \ar[r] \ar[d] &G\ar[d]\ar[r]&1\\
1\ar[r] & \pi_1^\top(\bar{S}',\bar{s}') \ar[r] & N_\phi\ar[r] & G\ar[r] &1
}
$$
The vertical arrows are surjective (left and middle: 
dense image, target discrete; right:
map is the identity).

\begin{Lemma}
$N_\phi$ is a fundamental group of a graph of groups.
\end{Lemma}

\begin{proof}
As $F:=\pi_1^\top(\bar{S}',\bar{s}')$ is free and of  finite 
index in $N_\phi$,
$\Omega'\subseteq\P^1_{\Cp}$ is the set of ordinary points of $N_\phi$. 
If $\Delta':=\Delta(\Omega')$ is the skeleton of $\Omega'$, then there is a
retraction map $\Omega'\to\Delta'$ \cite[4.1.6]{Ber90}, and $F$
acts on $\Delta'$ \cite[4.1.8]{Ber90}. $N_\phi$ does the same, and 
viewing $\Delta'$ as a graph, this leads to the following commuting 
diagramme
$$
\xymatrix{
\Omega' \ar[r] \ar[d] & \bar{S}' \ar[dr] \ar[d] & \\
\Delta'\ar[r]^{\!\!\text{univ}}\ar[drr]_{/N_\phi} &\Gamma' \ar[dr] &\cS\ar[d] \\
 &&\Gamma(N_\phi)
}
$$
The vertical arrows are all retraction maps, $\Gamma'$ is a graph with
trivial vertex and edge stabilisers, and $\Gamma(N_\phi)$ is the graph of
groups $(\Delta'/N_\phi,N_{\phi\,\bullet})$. Because $\Delta'$ is
a tree, $N_\phi\cong\pi^\BS(\Gamma(N_\phi))$.
\end{proof}

As we can replace $\Gamma(N_\phi)$ by the Kato graph $\Gamma^*(N_\phi)$
of (a representation of) $N_\phi$, we obtain

\begin{Lemma}	\label{lattice}
Let $\cS$ be a Mumford orbifold and 
$\phi\colon\bar{S}'\stackrel{/G}{\Pfeil}\cS$ a global chart with a Mumford
curve $\bar{S}'$ together with its universal topological cover
$\psi\colon\Omega'\to\bar{S}'$.
Then the following diagramme commutes and has exact rows and columns:
$$
\xymatrix{
&1\ar[d]&1\ar[d]&1\ar[d]&\\
1\ar[r]&\pi_1^\temp(\Omega,\bar{\omega})\ar[r]\ar[d]&\pi_1^\temp(\Omega,\bar{\omega})\ar[d]\ar[r]&1\ar[r]\ar[d]&1\\
1\ar[r]&\pi_1^\temp(\bar{S}',\bar{s}')\ar[d]\ar[r]&\pi_1^\orb(\mathcal{S},\bar{s})\ar[d]\ar[r]&G\ar[d]\ar[r]&1\\
1\ar[r]&\pi_1^\top(\bar{S}',\bar{s}')\ar[d]\ar[r]&
	\pi_1^\BS(\Gamma^*(N))\ar[d]\ar[r]&G\ar[d]\ar[r]&1\\
&1&1&1&
}
$$
Here, $N$ is the deck group of $\phi\circ\psi\colon\Omega\to\cS$.
\end{Lemma}

\begin{proof}
The two bottom rows are identical to the diagramme with exact rows
above, only $N_\phi$ has been replaced by the Bass-Serre
fundamental group $\pi_1^\BS(\Gamma^*(N_\phi))$. The middle  and left
columns are exact by the fundamental exact sequence in Lemma \ref{fundexact}.
The commutativity of the upper two boxes is clear, as the upper left
horizontal arrow is the identity map.
\end{proof}

Now fix a chart $\phi$ as above, and let
$$
\fT(\phi):=\left\{
			\raisebox{10mm}{
			\xymatrix{\pi_1^\orb(\mathcal{S},\bar{s}) \ar[r] \ar[d] & \PGL_2(\C_p)\\
			N_\phi \ar[ur]_{\text{faithful} \atop \text{discont.}}&
			}}
		    \right\}
$$
be the set of all representations of $\porbS$ factorising over faithfully
discontinuous representations of $N_\phi$. The isomorphism theorem shows
that the induced map
$$
\fT(\phi)\to\fT(N_\phi)
$$
is bijective, so $\fT(\phi)$ inherits from $\fT(N_\phi)$ the structure
of a (generally unconnected) manifold.

\smallskip
Let us fix the finite Galois group $G$, and let $\fK_G(\cS)$ be the set
of all isomorphism classes of global charts 
$\phi\colon\bar{S}'\stackrel{/G}{\Pfeil}\cS$ with $\bar{S}'$ a Mumford
curve and $\Aut(\bar{S}',\cS)\cong G$. According to the 
Riemann-Hurwitz formula, the genus $g'$ of $\bar{S}'$ is the same for
all $\phi\in\fK_G(\cS)$. If $g$ denotes the genus of $\bar{S}$, we 
can define:
$$
\fH_G(\cS):=\bigcup\limits_{\phi\in\fK_G(\cS)}\fM(N_\phi,F_{g'}),
$$ 
where $F_{g'}$ is the free group in $g'$ generators
sitting in the exact sequence
$$
1\to F_{g'}\to N_\phi\to G\to 1
$$
The Hurwitz space for covers of Mumford orbifolds of genus $g$ is obtained
in the following way: let
$$
\MOrb_g(\und{e})
$$
be the sub-category of the category of orbifolds consisting of the 
one-dimensional Mumford orbifolds of genus $g$ with ramification indexing
sequence $\und{e}:=(e_1,\dots,e_n)$.

\begin{Def}
The space 
$$
\fH_g(G;e_1,\dots,e_n):=\fH_g(G,\und{e}):=
\bigcup \limits_{\cS\in\MOrb_g(\und{e})} \fH_G(\cS)
$$
is called the {\em Mumford-Hurwitz space} for $g$-covers with 
{\em signature} $\und{e}$ of Mumford curves of genus $g$.
\end{Def}

In the following subsection, we will show that $\fK_G(\cS)$
is essentially finite, and that this gives only finitely many components
for the Mumford-Hurwitz space.

\subsection{The chart number}

Let $\Gamma^*:=\Gamma^*(N)$ be a Kato graph. 
A {\em stable model} of $\Gamma^*(N)$
is the graph obtained by stabilising the finite part of the Kato graph,
i.e.~by contracting edges $e$ in 
$\Gamma^*_\fin:=\Gamma^*\setminus\partial\Gamma^*$,
whenever $N_e$ is isomorphic to the stabiliser $N_v$ of an extremity $v$
of $e$ and  the valency of $v$ in $\Gamma^*_\fin$ is less than $3$.

\bigskip
We will call, for convenience, the number of stable
models for Kato graphs of charts $\phi\in\fK_G(\cS)$ the {\em chart number}
for the Mumford orbifold 
$$
\cS=(\bar{S},(\zeta_1,e_1),\dots,(\zeta_n,e_n))
$$ 
with Galois group $G$ and $n$ marked points.

\bigskip
Quite obvious, but essential, is 

\begin{Lemma}	\label{stabbound}
If $\phi\colon \bar{S}'\stackrel{/G}{\Pfeil}\cS$ is a finite global chart
for $\cS$, and if $\bar{S}'$ is a Mumford curve, then all vertex (and edge)
stabilisers of the corresponding Kato graph are subgroups of $G$.
\end{Lemma}

\begin{proof}
With the notations from the previous subsection, the diagramme
$$
\xymatrix{
T \ar[r]^{/F} \ar[dr]_{/N_\phi} & \Gamma \ar[d]^{/G}\\
& \Gamma^*(N_\phi)
}
$$
shows that vertex stabilisers occur only in the vertical map,
because $F$ is a free group: The maps with source $T$ are
both universal covers of the graph without groups $\Gamma$ on the one hand,
and of the Kato graph on the other hand.  
\end{proof}

\begin{Thm}	\label{finitechart}
The total  chart number of all $\cS\in\MOrb_g(\und{e})$ is finite,
if stable models of Kato graphs are counted without multiplicities.
\end{Thm}

\begin{proof}
We shall show that the number of vertices for  graphs of fixed Betti number
with $n$ cusps is
bounded from above. As the size of vertex stabilisers is also bounded by Lemma
\ref{stabbound}, this proves the theorem.

\medskip
Let first $\cS$ be a rational Mumford orbifold, $\phi\colon\bar{S'}\to\cS$
a chart with deck group $G$, and $\bar{S}'$ a Mumford curve.
Abbreviating $N:=N_\phi$,
we assume further that the Kato tree $\Gamma^*:=\Gamma^*(N)$
is  stable
and that $n\ge 3$.

\smallskip
In case $\Gamma^*$ is irreducible, i.e.~no edge stabiliser is
trivial, \cite[Proposition 15]{Brad02a} implies that the valency of
any vertex in $\Gamma^*$ is at most $3$. 

\smallskip
Let $v$ be a vertex in $\Gamma^*_\fin$, and $\val(v)$ its valency in
that tree. There are the following possibilities:
\begin{enumerate}
\item $\val(v)=1$. 
\item $\val(v)=2$ and the $N_{e}$ for both edges  $e$ emanating from
$v$ are cyclic.
\item $\val(v)=2$ and $N_{e}$ for an edge $e$ emanating from $v$
is non-cyclic.
\item $\val(v)=3$.
\end{enumerate}

For proving boundedness in each case, we rely on how
$\Gamma^*$ is obtained by glueing Kato trees for finite groups.
This is explained in the proof of \cite[Theorem 2]{Brad02a}.
These trees have two cusps if the group is cyclic, and three cusps otherwise.

\smallskip
1. In this case, $v$ has a cusp emanating from it in $\Gamma^*$. 

\smallskip
2. Here, $N_{v}$ is not cyclic, and $v$ has a cusp in $\Gamma^*$. 

\smallskip
3. Let $v=o(e)$, where $e$ is an edge with $N_{e}$ not cyclic.
Then $v(t(e))=3$ in $\Gamma^*$, because $N_{t(e)}$
is not cyclic. So, here is a cusp going out of $v$.

\smallskip
4. For each vertex $v$ with $\val(v)=3$, there is at least one
extremal vertex in $\Gamma^*_\fin$. For such vertices, we are in case 1. 

\smallskip
From this, we see that the number of vertices in $\Gamma^*$ is bounded.

\medskip
In case $\Gamma^*$ contains edges with trivial stabiliser,
each maximal subtree without trivial edge groups has at least one cusp
going out. Such a subtree is called an {\em irreducible component}
of $\Gamma^*$. The number of irreducible components is therefore bounded.

\medskip
Let now $\cS$ be of arbitrary genus.
The considerations above for a fundamental domain, viewed as a tree
of groups, of $\Gamma^*$
in its universal covering tree prove the theorem in this general case.
\end{proof}


\begin{Cor}	\label{MH}
The Mumford-Hurwitz space $\fH_g(G,\und{e})$ is a coarse moduli space
parametrising $G$-covers with signature $\und{e}=(e_1,\dots, e_n)$ 
of Mumford covers of 
genus $g$, and has only finitely many
components of equal dimension $3g-3+n$.   
\end{Cor}

\begin{proof}
The moduli space property is given by Proposition \ref{moduli}.
Equidimensionality has been proven in Proposition \ref{equidim}.
The  
finiteness property of the Mumford-Hurwitz space
follows from  the 
connectedness
of each $\fM(N_\phi)$ (of which $\fM(N_\phi, F_{g'})$ is a finite cover)
and from the finiteness of the chart number (Theorem \ref{finitechart}).
\end{proof}

\subsection{Examples}

\subsubsection{Tate orbifolds}	\label{Tate}

Here we study Galois covers of the Tate orbifold 
$\cS_\lambda=(\P_1,(0,2), (1,2),(\infty,2),(\lambda,2))$. This means that 
$G=C_2$, the cyclic group of order two. Let $p\neq 2$. The rational orbifold
$\cS_\lambda$ can be defined for all $\lambda\neq 0,1,\infty$, but
it is a global orbifold chart of Mumford type
 if and only if $|\lambda||\lambda-1|\neq 1$. In fact, the equation
$$
y^2=x(x-1)(x-\lambda)
$$
defines a chart $\phi\colon T\to\cS_\lambda$ with Galois group $C_2$,
and $T$ is known to be a Tate elliptic curve if and only if
$|\lambda||\lambda|\neq 1$.
We have
$$
\fH_0(C_2;2,2,2,2)
=\{\lambda\in\Cp^\times\setminus\{1\}\mid|\lambda'||\lambda'-1|\neq 1\}
=\fM(C_2*C_2,\Z).
$$
The moduli space $\fM(C_2*C_2,\Z)$ is given by the chart
$$
\Cp^\times\stackrel{/C_2*C_2}{\Pfeil}\cS_\lambda
$$
coming from the universal topological cover of the Tate curve;
it is a Galois cover whose group $\pi_1^\BS(\Gamma)\cong C_2*C_2$ 
comes from the  cover of analytic skeletons viewed as graph with groups
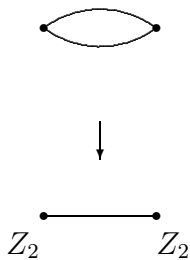
\begin{center}
\setlength{\unitlength}{.5cm}
\begin{picture}(5,8)
\put(1,2){\circle*{.2}}
\put(4,2){\circle*{.2}}
\put(1,2){\line(1,0){3}}
\put(1,7){\circle*{.2}}
\put(4,7){\circle*{.2}}
\put(2.5,4.5){\vector(0,-1){1}}
\put(2.5,5){\arc(1.5,2){71}}
\put(2.5,9){\arc(-1.5,-2){71}}
\put(0,1){$Z_2$}
\put(4,1){$Z_2$}
\end{picture}
\end{center} 

In order to understand the associated surjection 
$\pi_1^\orb(\cS_\lambda,\bar{s})\to\pi_1^\BS(\Gamma)$, we observe

\begin{Prop}	\label{invertgen}
Let $\cS_\lambda$ be as above. Then:
\begin{enumerate}
\item
If  
$\cS_\lambda$ is a Tate orbifold, then
the isomorphisms 
$$
\pi_1^\orb(\cS_\lambda,\bar{s})\cong (\Z\times\hat\Z)\rtimes C_2
\cong\hat\Z\rtimes(C_2*C_2)
$$
hold.
\item
If $\cS_\lambda$ is not a Tate orbifold, then 
$$
\pi_1^\orb(\cS_\lambda,\bar{s})\cong \hat\Z^2\rtimes C_2.
$$
\end{enumerate}
In both cases, $C_2$ inverts 
each of the two topological generators of $\Z\times\hat\Z$, resp.~$\hat\Z^2$,
and in the first case
 both generators of $C_2*C_2$ in the group on the right
invert the topological generator
of $\hat\Z$. 
\end{Prop}

\begin{proof}
The first isomorphy in 1.~as well as the isomorphy of 2.~is shown in \cite[Remarks 4.5.6.(c)]{AndIII}, while the 
second isomorphy in 1.~follows from  Lemma \ref{lattice}: the diagramme is in our case
$$
\xymatrix{
&1\ar[d]&1\ar[d]&1\ar[d]&\\
1\ar[r]&\hat{\Z}\ar[r]\ar[d]&\hat{\Z}\ar[d]\ar[r]&1\ar[r]\ar[d]&1\\
1\ar[r]&\Z\times\hat{\Z}\ar[d]\ar[r]&\pi_1^\orb(\mathcal{S}_\lambda,\bar{s})\ar[d]\ar[r]&Z_2\ar[d]\ar[r]&1\\
1\ar[r]&\Z\ar[d]\ar[r]&Z_2\ast Z_2\ar[d]\ar[r]&Z_2\ar[d]\ar[r]&1\\
&1&1&1&
}
$$ 
For the action of $C_2$, let $\phi_n\colon \Omega_n\to T$ be a tempered
cover of the elliptic curve $T$. 
 In the case that $T$ has good reduction, $\phi_n$ is an isogeny.
If $T$ is a Tate curve, we may assume that $\phi_n$ is the composition
of the universal cover of a Tate curve with an isogeny.
In either case, $\phi_n$ is a homomorphism of Abelian groups.
Let $\sigma_n$ be the involution of $\Omega_n$ lifting the involution 
$\sigma\colon t\to t^{-1}$ of $T$ (we shall write $T$ as a multiplicative 
Abelian group for obvious reasons).
Any other lift of $\sigma$ is of the form
$$
\sigma_n\cdot\epsilon_n\colon\Omega_n\to\Omega_n,
\;\omega\mapsto \omega^{-1}\cdot\epsilon_n(\omega)
$$
with $\epsilon_n(\omega)\in\ker\phi_n$.
As $\epsilon_n$ is continuous and $\ker\phi_n$ discrete, we find
$$
\epsilon_n=\const=\sigma_n(1).
$$
Now, 
$\sigma:=(\sigma_n)\in
\pi_1^\orb(\cS_\lambda,\bar{s})\setminus\pi_1^\temp(T,\bar{1})$,
$\epsilon:=(\epsilon_n)\in\pi_1^\temp(T,\bar{1})$,
and the involution $\sigma$ acts on the Abelian group $\pi_1^\temp(T,\bar{1})$
by conjugation
$$
\pi_1^\temp(T,\bar{1})\to\pi_1^\temp(T,\bar{1}),\;
\gamma\mapsto(\sigma\epsilon)\cdot\gamma\cdot(\sigma\epsilon)^{-1}
=\sigma\gamma^{-1}\stackrel{!}{=}\gamma^{-1}.
$$
The last equality can be checked for each $\phi_n$: Here the fibre
$\phi_n^{-1}(1)$ equals the Abelian deck group $G_n$. Then each $\gamma_n$ from
$\gamma=(\gamma_n)\in\pi_1^\temp(T,\bar{1})$ is the translation 
$x\mapsto x\cdot\gamma_n$, and for all $x\in G_n$ we have
$$
\gamma_n^{\sigma_n}(x)=\sigma_n\gamma_n\sigma_n(x)=\sigma_n\gamma_n(x^{-1})=\sigma_n(x^{-1}\cdot\gamma_n)=\gamma_n^{-1}\cdot x=\gamma_n^{-1}(x)\,,
$$
in other words, $\gamma^{\sigma_n}=\gamma_n^{-1}$.

\smallskip
As $C_2*C_2$ is the free product of two copies of $\langle\sigma\rangle$,
the exact diagramme from the beginning of the proof gives the last assertion. 
\end{proof}

\begin{Rem}
If $T$ has good reduction, then the temperate fundamental group
is the algebraic fundamental group, 
$\pi_1^\temp(T,\bar{1})\cong\pi_1^\alg(T,\bar{1})$, as temperate covers
of $T$ are all finite. In this case, it is known that $C_2$ 
inverts each of the two topological generators of $\hat\Z^2$, for example
by using the fact, that this is the case in the complex situation
for the (complex) topological fundamental group and
 that the pro-finite completion of the latter is
the algebraic fundamental group.
Thus, Proposition {\rm\ref{invertgen}} can be viewed as 
a generalisation of the good
reduction case.
\end{Rem}

\subsubsection{Triangle groups}

A special type of rational Mumford orbifolds are 
{\em Mumford-Schwarz orbifolds} 
$$
\cS(e_0,e_1,e_\infty)=(\P^1_{\Cp},(0,e_0),(1,e_1),(\infty,e_\infty)).
$$ 
Such an orbifold $\cS(e_0,e_1,e_\infty)$ is a quotient
of a Mumford curve $X$, and if the genus of $X$ is $\ge 2$, it 
 is called {\em hyperbolic}. The Bass-Serre
fundamental group of a corresponding Kato graph is called a 
{\em $p$-adic triangle group of Mumford type} in the hyperbolic case.   

\medskip
Kato has proven that $p$-adic triangle groups of Mumford type exist only
for $p\le 5$ \cite{Kat00}, more exactly:

\begin{Thm}[Kato, 1999]	\label{Kato}
There are infinitely many hyperbolic Mumford-Schwarz orbifolds
for $p\le 5$, and none for $p>5$.
For given ramification indices $e_0$, $e_1$ 
and $e_\infty$ there are only finitely many $p$-adic triangle
groups $\Delta(e_0,e_1,e_\infty)$ of Mumford type. 
\end{Thm}

This also implies the finiteness of the chart number for Mumford-Schwarz
orbifolds, 
a special case of Theorem \ref{finitechart}.
In any case, we partially recover from Corollary \ref{MH}:

\begin{Cor}
The Mumford-Hurwitz space $\fH_0(G; e_0, e_1, e_\infty)$
is a finite set. 
\begin{enumerate}
\item
If $p>5$, then it is non-empty, only if $G$ equals 
one of the finite $\Delta(e_0,e_1,e_\infty)$. 
\item
It is non-empty for finite quotients of infinitely many
infinite triangle groups, if $p\le 5$.
\end{enumerate}
\end{Cor}

\subsubsection{Quadrangle groups}

The formula in \cite[Theorem 2]{Brad02a} and its proof imply that the
only way to obtain Kato tree with four cusps is by pasting two
three-cusped Kato graphs along an elementary Kato graph with two cusps:
$$
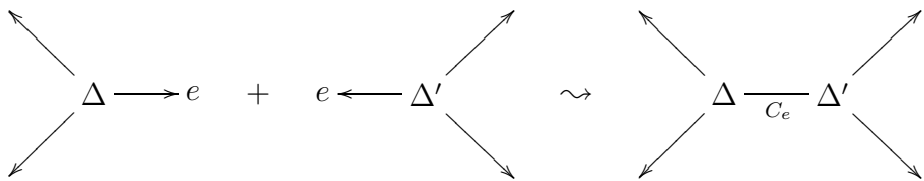
\xymatrix{
 & & \\
 &\Delta \ar[ul]\ar[dl] \ar[r]& e\\
 & &
}
\raisebox{-1.2cm}{\hspace{5mm}$+$}
\hspace{5mm}
\xymatrix{
 &&\\
e & \Delta' \ar[ur]\ar[dr]\ar[l]\\
&&
}
\raisebox{-1.2cm}{$\quad\leadsto\quad$}
\xymatrix{
 &&&\\
&\Delta\ar[ul]\ar[dl]\ar@{-}[r]_{C_e} & \Delta' \ar[ur]\ar[dr]&\\
&&&
}
$$
By Theorem \ref{Kato}, $\Delta$ and $\Delta'$ are both finite, if $p>5$.
In other words, The fundamental group is an amalgam of two
finite groups. This has also been observed in
\cite{vdPV01} by different means.

\smallskip
In the case that $p\le 5$, one has only to check in Kato's list of amalgams
\cite[Theorem (2)]{Kat00}, where exactly the cusps emanate and what
their decomposition groups are (an easy task). 
Then ``glue together'' two cusps
with equal stabiliser $C_e$ (as then one has an allowed segment of
Herrlich's list \cite{Her82}).

\subsubsection{Cyclic covers}

From Lemma \ref{stabbound} we get immediately

\begin{Prop}	\label{cyclic}
Let $\cS$ be any Mumford orbifold and $\phi\in\fK_{C_n}(\cS)$. Then
all vertex groups in the Kato graph for $\phi$ are cyclic, and all
edge groups are either trivial or equal to the stabiliser of 
one of its extremities.
\end{Prop}

The Kato graph reveals the relative position of the branch points
\cite[Theorem 4]{Brad02a},
so we have

\begin{Cor}
If a Mumford orbifold $\cS=(\bar{S},(\zeta_1,e_1),\dots,(\zeta_m,e_m))$
admits a chart $\phi\in\fK_{C_n}(\cS)$, then $m$ is even,
and the branch points are separated by a pure affinoid covering of $\cS$
into pairs of branch points
with equal decomposition group
\end{Cor}

The meaning of Proposition \ref{cyclic} in case $\cS$ is of genus zero 
is that, 
if one stabilises 
the irreducible
components of a Kato tree of $\cS$, there remains only one vertex,
i.e.~one has an elementary Kato tree for the corresponding decomposition group.

\medskip
In order to find the relative distances of the pairs of branch points in
$\cS$, one has to calculate the lengths of the edges
joining the vertices from which cusps emanate. The subtleties lie in
the case that $p$ divides the order of a vertex group.
This dealt with in Sections 6.3 and 6.4 of \cite{Brad02}.

\subsection*{Acknowledgements}

The author expresses his gratitude to Y.~Andr\'e for sending him
an early version of \cite{AndIII}, and to F.~Kato for inventing
his graphs. He thanks S.~K\"uhnlein for valuable discussions
on temperate covers of Tate curves and his thesis advisor
F.~Herrlich for patiently answering numerous  questions.

{\sc Universit\"at Karlsruhe, Mathematisches Institut I, Englerstr.~2, D-76128
Karlsruhe, Germany}

e-mail: {\tt bradley@math.uni-karlsruhe.de}



\begin{thebibliography}{ABCDE}


\bibitem[AndIII]{AndIII}Yves Andr\'e. {\em $p$-adic orbifolds and monodromy}, chapter III from {\em Period mappings and differential equations. From $\C$ to $\C_p$}, T\^{o}hoku-Hokkaid\^{o} Lectures in Arithmetic Geometry, Preprint {\tt math.AG/0203194}

\bibitem[Ber90]{Ber90}Vladimir G.~Berkovich. {\em Spectral Theory and Analytic Geometry over Non-Archimedean Fields}, Mathematical Surveys and Monographs, Number 33, AMS (1990)



\bibitem[Brad02]{Brad02}Patrick Erik Bradley. {$p$-adische Hurwitzr\"aume}, dissertation, Universit\"at Karlsruhe (2002)

\bibitem[Brad02a]{Brad02a}Patrick Erik Bradley. {On Mumford Orbifolds}, preprint {\tt math.AG/0207299}



\bibitem[dJ95]{dJ95}Aise Johan de Jong. {\em \'Etale fundamental groups of non-archimedean analytic spaces}, Compositio Math.\ 97, 89-118 (1995)



\bibitem[Ger81]{Ger81}Lothar Gerritzen. {\em Zur analytischen Beschreibung des Raumes der Schottky-Mumford-Kurven}, Math.\ Ann.\ 255, 259-271 (1981)


\bibitem[GvP80]{GvP80}Lothar Gerritzen, Marius van der Put. {\em Schottky Groups and Mumford Curves}, Lecture Notes in Mathematics 817, Springer-Verlag, 1980



\bibitem[Her82]{Her82}Frank Herrlich. {\em $p$-adisch diskontinuierlich einbettbare Graphen von Gruppen}, Arch.\ Math., Vol.\ 39, 204-216 (1982)

\bibitem[Her84]{Her84}Frank Herrlich. {\em On the stratification of the moduli space of mumford curves}, Groupes d'\'etude d'Analyse ultram\'etrique (Y.\ Amice, G.\ Christol, P.\ Robba) 11e ann\'ee, 1983/84, no.\ 18

\bibitem[Her85]{Her85}Frank Herrlich. {\em Nichtarchimedische Teichm\"ullerr\"aume}, Habilitationsschrift  Bochum (1985)

\bibitem[Her87]{Her87}Frank Herrlich. {\em Nichtarchimedische Teichm\"ullerr\"aume}, Indag.\ Math.\ 49, 145-169 (1987)

\bibitem[Kat00]{Kat00}Fumiharu Kato. {\em $p$-adic Schwarzian triangle groups of Mumford type}, preprint {\tt math.AG/9908174}

\bibitem[Kat01]{Kat01}Fumiharu Kato. {\em Non-archimedean orbifolds covered by Mumford curves}, preprint (2001) 


\bibitem[Khr91]{Khr91}D.\ G.\ Khramtsov. {\em Finite graphs of groups with isomorphic fundamental groups}, Algebra Logic 30, No.5, 389-409 (1991); translation from Algebra Logika 30, No.5, 595-623 (1991)


\bibitem[L\"ut83]{Lut83}Werner L\"utkebohmert. {\em Ein globaler Starrheitssatz f\"ur Mumfordkurven}, J.\ f\"ur reine und angew.\ Math.\ 340, 118-139 (1983)















\bibitem[vdPV01]{vdPV01}Marius van der Put, Harm Voskuil. {\em Discontinuous subgroups of $\PGL_2(K)$}, preprint so far available only on van der Put's Homepage {\tt http://www.math.rug.nl/\~{}vdput/ } (2001)


\end{thebibliography}
\end{document}